\documentclass[a4paper,11pt]{amsart}
\usepackage{amsmath,amssymb,amscd,graphicx,psfrag,epsf,amsthm,times,hyperref}
\usepackage[all]{xy}

\DeclareFontFamily{OMS}{rsfs}{\skewchar\font'60}
\DeclareFontShape{OMS}{rsfs}{m}{n}{<-5>rsfs5 <5-7>rsfs7 <7->rsfs10 }{}
\DeclareSymbolFont{rsfs}{OMS}{rsfs}{m}{n}
\DeclareSymbolFontAlphabet{\scr}{rsfs}

\newcommand{\sC}{\scr{C}}
\newcommand{\sD}{\scr{D}}
\newcommand{\sE}{\scr{E}}
\newcommand{\sF}{\scr{F}}
\newcommand{\sG}{\scr{G}}
\newcommand{\sH}{\scr{H}}

\newcommand{\sL}{\scr{L}}
\newcommand{\sM}{\scr{M}}
\newcommand{\sN}{\scr{N}}
\newcommand{\sO}{\scr{O}}

\newcommand{\sT}{\scr{T}}

\newcommand{\bN}{\mathbb{N}}

\newcommand{\bP}{\mathbb{P}}

\newcommand{\bZ}{\mathbb{Z}}

\newcommand{\cC}{\mathcal{C}}

\newcommand{\cH}{\mathcal{H}}

\newcommand{\cP}{\mathcal{P}}
\newcommand{\cQ}{\mathcal{Q}}

\newcommand{\cU}{\mathcal{U}}

\newcommand{\into}{\hookrightarrow}

\newcommand{\leteq}{\colon\!\!\!=}

\DeclareMathOperator{\codim}{codim}

\DeclareMathOperator{\im}{{im}}

\DeclareMathOperator{\Proj}{{Proj}}

\DeclareMathOperator{\rk}{{rk}}

\DeclareMathOperator{\supp}{{supp}}

\DeclareMathOperator{\Sym}{{Sym}}

\def\coh#1.#2.#3.{H^{#1}(#2,#3)}
\def\dimcoh#1.#2.#3.{h^{#1}(#2,#3)}
\newcommand{\union}\cup
\newcommand{\intersect}\cap
\newcommand{\Union}\bigcup
\newcommand{\Intersect}\bigcap
\newcommand{\assign}\longmapsto
\newcommand{\resto}{\big\vert_}

\def\factor#1.#2.{\left. \raise 2pt\hbox{$#1$} \right/
  \hskip -2pt\raise -2pt\hbox{$#2$}}

\newcommand\mypagesizea{
\textwidth= 6.5in
\textheight=9.5in
\voffset-.15in
\hoffset-.75in
\marginparwidth=56pt
}

\newtheoremstyle{bozont}{3pt}{3pt}%
     {\itshape}
     {}
     {\bfseries}
     {.}
     {.5em}
     {\thmname{#1}\thmnumber{ #2}\thmnote{ \rm #3}}
\newtheoremstyle{bozont-sf}{3pt}{3pt}%
     {\itshape}
     {}
     {\sffamily}
     {.}
     {.5em}
     {\thmname{#1}\thmnumber{ #2}\thmnote{ \rm #3}}
\newtheoremstyle{bozont-sc}{3pt}{3pt}%
     {\itshape}
     {}
     {\scshape}
     {.}
     {.5em}
     {\thmname{#1}\thmnumber{ #2}\thmnote{ \rm #3}}
\newtheoremstyle{bozont-remark}{3pt}{3pt}%
     {}
     {}
     {\scshape}
     {.}
     {.5em}
     {\thmname{#1}\thmnumber{ #2}\thmnote{ \rm #3}}
\newtheoremstyle{bozont-def}{3pt}{3pt}%
     {}
     {}
     {\bfseries}
     {.}
     {.5em}
     {\thmname{#1}\thmnumber{ #2}\thmnote{ \rm #3}}
\newtheoremstyle{bozont-reverse}{3pt}{3pt}%
     {\itshape}
     {}
     {\bfseries}
     {.}
     {.5em}
     {\thmnumber{#2.}\thmname{ #1}\thmnote{ \rm #3}}
\newtheoremstyle{bozont-reverse-sc}{3pt}{3pt}%
     {\itshape}
     {}
     {\scshape}
     {.}
     {.5em}
     {\thmnumber{#2.}\thmname{ #1}\thmnote{ \rm #3}}
\newtheoremstyle{bozont-reverse-sf}{3pt}{3pt}%
     {\itshape}
     {}
     {\sffamily}
     {.}
     {.5em}
     {\thmnumber{#2.}\thmname{ #1}\thmnote{ \rm #3}}
\newtheoremstyle{bozont-remark-reverse}{3pt}{3pt}%
     {}
     {}
     {\sc}
     {.}
     {.5em}
     {\thmnumber{#2.}\thmname{ #1}\thmnote{ \rm #3}}
\newtheoremstyle{bozont-def-reverse}{3pt}{3pt}%
     {}
     {}
     {\bfseries}
     {.}
     {.5em}
     {\thmnumber{#2.}\thmname{ #1}\thmnote{ \rm #3}}
\newtheoremstyle{bozont-def-newnum-reverse}{3pt}{3pt}%
     {}
     {}
     {\bfseries}
     {}
     {.5em}
     {\thmnumber{#2.}\thmname{ #1}\thmnote{ \rm #3}}
\theoremstyle{bozont} 

\newtheorem{thm}{Theorem}

\newtheorem{lemma}[thm]{Lemma}
\newtheorem{cor}[thm]{Corollary}

\newtheorem{prop}[thm]{Proposition}

\newtheorem*{thm*}{Theorem}

\theoremstyle{definition}
\newtheorem{defn}[thm]{Definition}

\newtheorem{ack}{Acknowledgments}

\newtheorem{defn-thm}[thm]{Definition-Theorem}  

\newtheorem{rem}[thm]{Remark}

\theoremstyle{remark}

\newtheorem*{not-and-def}{Notation and definitions}

\newcommand{\p}[0]{{\mathbb P}}

\newcommand{\z}[0]{{\mathbb Z}}

\newcommand{\m}[0]{{\mathfrak m}}

\newcommand{\q}[0]{{\mathbb Q}}
\newcommand{\map}[0]{\dasharrow}

\newcommand{\rat}[0]{\operatorname{RatCurves}^n}

\newcommand{\chow}[0]{\operatorname{Chow}}

\newcommand{\rup}[1]{\lceil{#1}\rceil}

\def\into{\DOTSB\lhook\joinrel\rightarrow}

\numberwithin{thm}{section}
\numberwithin{equation}{thm}

\newcommand\dual{\ast}

\mypagesizea

\begin{document}

\title[Characterizations of projective spaces and hyperquadrics]{Cohomological
  characterizations of projective spaces and hyperquadrics}

\author{Carolina Araujo, St\'ephane Druel, S\'andor J.\ Kov\'acs} 

\address{Carolina Araujo: \sf IMPA, Estrada Dona Castorina 110, Rio de
  Janeiro, 22460-320, Brazil} 
\email{caraujo@impa.br}

\address{\vskip -.5cm St\'ephane Druel: \sf Institut Fourier, UMR 5582 du
  CNRS, Universit\'e Joseph Fourier, BP 74, 38402 Saint Martin
  d'H\`eres, France} 
\email{druel@ujf-grenoble.fr}

\address{\vskip -.5cm S\'andor J.\ Kov\'acs: \sf University of
  Washington, Department of Mathematics, Box 354350, Seattle, WA
  98195, USA} 
\email{kovacs@math.washington.edu}

\date{July 29, 2007} 

\thanks{The first named author was partially supported by a CNPq Research 
  Fellowship.}

\thanks{The second named author was partially supported by the 3AGC
  project of the A.N.R}

\thanks{The third named author was supported in part by NSF Grant
  DMS-0554697 and the Craig McKibben and Sarah Merner Endowed
  Professorship in Mathematics.}


\maketitle




\section{Introduction}
\label{intro}

Projective spaces and hyperquadrics are the simplest projective
algebraic varieties, and they can be characterized in many ways.  The
aim of this paper is to provide a new characterization of them in
terms of positivity properties of the tangent bundle
(Theorem~\ref{main_thm}).

The first result in this direction was Mori's proof of the Hartshorne
conjecture in \cite{mori79} (see also Siu and Yau \cite{siu-yau}),
that characterizes projective spaces as the only manifolds having
ample tangent bundle.  Then, in \cite{wahl}, Wahl characterized
projective spaces as the only manifolds whose tangent bundles contain
ample invertible subsheaves.  Interpolating Mori's and Wahl's results,
Andreatta and Wi\'sniewski gave the following characterization:

\begin{thm*}[\cite{andreatta_wisniewski}]
  Let $X$ be a smooth complex projective $n$-dimensional variety.
  Assume that the tangent bundle $T_X$ contains an ample locally free
  subsheaf $\sE$ of rank $r$.  Then $X\simeq \bP^n$ and either
  $\sE\simeq \sO_{\bP^n}(1)^{\oplus r}$ or $r=n$ and $\sE= T_{\bP^n}$.
\end{thm*}

We note that earlier, in \cite{campana_peternell}, Campana and
Peternell obtained the same result for $r\geq n-2$.

Let $\sE$ be an ample locally free subsheaf of $T_{\bP^n}$ of rank
$p<n$.  By taking its determinant, we obtain a non-zero section in
$H^0(\bP^n,\wedge^{p} T_{\bP^n}\otimes \sO_{\bP^n}(-p))$.  On the other
hand, most sections in $H^0(\bP^n,\wedge^{p} T_{\bP^n}\otimes
\sO_{\bP^n}(-p))$ do not come from ample locally free subsheaves of
$T_{\bP^n}$.

This motivates the following characterization of
projective spaces and hyperquadrics, which was conjectured by
Beauville in \cite{beauville_symplectic}. Here $Q_p$ denotes a smooth
quadric hypersurface in $\bP^{p+1}$, and $\sO_{Q_{p}}(1)$ denotes the 
restriction of $\sO_{\bP^{p+1}}(1)$ to $Q_p$.
When $p=1$, $(Q_{1},\sO_{Q_{1}}(1))$ is just $(\bP^{1},\sO_{\bP^{1}}(2))$.

\begin{thm}\label{main_thm} 
  Let $X$ be a smooth complex projective $n$-dimensional variety and
  $\sL$ an ample line bundle on $X$.  If $H^0(X,\wedge^{p}
  T_{X}\otimes \sL^{-p})\neq 0$ for some positive integer $p$, then
  either $(X,\sL)\simeq (\bP^{n},\sO_{\bP^{n}}(1))$, or $p=n$ and
  $(X,\sL)\simeq(Q_{p},\sO_{Q_{p}}(1))$.
\end{thm}

The statement of this theorem can be interpreted in the following way.
Let $X$ be a smooth complex projective $n$-dimensional variety and
$\sL$ an ample line bundle on $X$.  Consider the sheaf $\sT_\sL\leteq
T_X\otimes \sL^{-1}$. Then Wahl's theorem \cite{wahl} says that if
$\coh 0. X. \sT_\sL.\neq 0$ then $X\simeq \bP^n$.
Theorem~\ref{main_thm} generalizes this statement to the case when one
only assumes that $\coh 0. X. \wedge^p\sT_\sL.\neq 0$ for some
$0<p\leq n$.

In order to prove Theorem~\ref{main_thm}, first notice that $X$ is
uniruled by \cite[Corollary 8.6]{miyaoka}.  Next observe that if the
Picard number of $X$ is $1$, then it is necessarily a Fano variety. If
the Picard number is larger than $1$, then we fix a minimal covering
family $H$ of rational curves on $X$, and follow the strategy in
\cite{andreatta_wisniewski} of looking at the $H$-rationally connected quotient
$\pi:X^\circ\to Y^\circ$ of $X$ (see
Section~\ref{section_minimal_rat_curves} for definitions).  We show
that any non-zero section $s\in H^0(X,\wedge^{p} T_{X}\otimes
\sL^{-p})$ restricts to a non-zero section $s^\circ\in
H^0(X^\circ,\wedge^{p} T_{X^\circ/Y^\circ}\otimes \sL^{-p})$, except
in the very special case when $p=2$ and $X\simeq Q_2$. This is
achieved in Section~\ref{section_fibration}. Afterwards we need to
deal with two cases: the case when $X$ is a Fano manifold with Picard
number $1$, and the case in which the $H$-rationally connected quotient
$\pi:X^\circ\to Y^\circ$ is either a projective space bundle or a
quadric bundle, and $H^0(X^\circ,\wedge^{p} T_{X^\circ/Y^\circ}\otimes
\sL^{-p})\neq 0$.

When $X$ is a Fano manifold with Picard number $\rho(X)=1$, the result
follows from the following.

\begin{thm}[(=\,{Theorem~\ref{beauville_for_rho=1}})]
  \label{thm:beauville_for_rho=1}
  Let $X$ be a smooth $n$-dimensional complex projective variety with
  $\rho(X)=1$, $\sL$ an ample line bundle on $X$, and $p$ a positive
  integer.  If $H^0(X,T_{X}^{\otimes p}\otimes \sL^{-p})\neq 0$, then
  either $(X,\sL)\simeq (\p^{n},\sO_{\p^{n}}(1))$, or $p=n\geq
  3$ and $(X,\sL)\simeq(Q_{p},\sO_{Q_{p}}(1))$.
\end{thm}


The paper is organized as follows.  In
Section~\ref{section_minimal_rat_curves} we gather old and new results
about minimal covering families of rational curves and their
rationally connected quotients.  In Section~\ref{section_relative_-K}
we show that the relative anti\-canonical bundle of a generically
smooth surjective morphism from a normal projective $\q$-Gorenstein
variety onto a smooth curve is never ample.  This will be used to
treat the case when the $H$-rationally connected quotient
$\pi:X^\circ\to Y^\circ$ is a quadric bundle.  In
Section~\ref{section_lifting_p_derivations}, we show that
$p$-derivations can be lifted to the normalization.  This technical
result will be used in the following section, which is the technical
core of the paper.  In Section~\ref{section_fibration}, we study the
behavior of non-zero global sections of bundles of the form
$\wedge^{p}T_{X}\otimes \sM$ with respect to fibrations $X\to Y$.  We
also prove some general vanishing results, such as the following.

\begin{thm}[(=\,{Corollary~\ref{vanishing}})]\label{thm:vanishing}
  Let $X$ be a smooth complex projective variety and $\sL$ an ample line bundle
  on $X$.  If $H^0(X,\wedge^{p}T_{X}\otimes \sL^{-p-1-k})\neq 0$ for
  integers $p\ge 1$ and $k\ge 0$, then $k=0$ and $(X,\sL)\simeq
  (\bP^{p},\sO_{\bP^{p}}(1))$.
\end{thm}
Finally, in Section~\ref{section_proof} we prove
Theorem~\ref{thm:beauville_for_rho=1} and put things together to prove
Theorem~\ref{main_thm}.

\begin{not-and-def}
  Throughout the present article we work over the field of complex
  numbers unless otherwise noted.  By a vector bundle we mean a locally free sheaf and by a
  line bundle an invertible sheaf. If $X$ is a variety and $x\in X$,
  then $\kappa(x)$ denotes the residue field
  $\factor{\sO_{X,x}}.{\mathfrak m_{X,x}}.$.
  Given a variety $X$, we denote by $\rho(X)$ the Picard number of
  $X$.  
  If $\sE$ is
  a vector bundle over a variety $X$, we denote by $\sE^*$ its dual
  vector bundle, and by $\bP(\sE)$ the Grothendieck projectivization
  $\Proj_X(\Sym (\sE))$. 
  For a morphism $f:X\to T$, the fiber of $f$ over $t\in T$ is denoted
  by $X_t$.
\end{not-and-def}

\begin{ack}
  The work on this project benefitted from support from various institutions and from
  discussions with some of our colleagues. In particular, the second and third named
  authors' visits to the \emph{Instituto Nacional de Matem\'atica Pura e Aplicada}
  and the first and second named authors' visit to the \emph{Korea Institute of
    Advanced Studies} were essentially helpful. The former visits were made possible
  by support from the {ANR}, the {NSF} and {IMPA}. The latter took place during a
  workshop organized by {Jun-Muk Hwang} with support from {KIAS}.  We would like to
  thank these institutions for their support and \emph{Jun-Muk Hwang} for his
  hospitality.  We would also like to thank \emph{J\'anos Koll\'ar} for helpful
  discussions and suggestions that improved both the content and the presentation of
  this article.
\end{ack}



\section{Minimal rational curves on uniruled varieties}

\label{section_minimal_rat_curves}

In this section we gather some properties of minimal covering families
of rational curves and their corresponding rationally connected
quotients.  For more details 
see \cite{kollar}, \cite{debarre}, or \cite{araujo_kollar}.

Let $X$ be a smooth complex projective uniruled variety and $H$ an
irreducible component of $\rat(X)$. Recall that only general points in
$H$ are in 1:1-correpondence with the associated curves in $X$.

We say that $H$ is a \emph{covering family of rational curves on $X$}
if the corresponding universal family dominates $X$.  A covering
family $H$ of rational curves on $X$ is called \emph{unsplit} if it is
proper.  It is called \emph{minimal} if, for a general point $x\in X$,
the subfamily of $H$ parametrizing curves through $x$ is proper.  As
$X$ is uniruled, a minimal covering family of rational curves on $X$
always exists.  One can take, for instance, among all covering
families of rational curves on $X$ one whose members have minimal
degree with respect to a fixed ample line bundle.

Fix a minimal covering family $H$ of rational curves on $X$.  Let $C$
be a rational curve corresponding to a general point in $H$, with
normalization morphism $f:\bP^1\to C\subset X$. We denote by $[C]$ or
$[f]$ the point in $H$ corresponding to $C$.  We denote by $f^*T_X^+$
the subbundle of $f^*T_X$ defined by
$$
f^*T_X^+ = \im \left[H^0\big(\bP^1, f^*T_X(-1)\big)\otimes
  \sO_{\bP^1}(1) \to f^*T_X\right] \into f^*T_X.
$$
By \cite[IV.2.9]{kollar}, if $[f]$ is a general member of $H$, then
$f^*T_X\simeq \sO_{\bP^1}(2)\oplus \sO_{\bP^1}(1)^{\oplus d}\oplus
\sO_{\bP^1}^{\oplus (n-d-1)}$, where $d=\deg(f^*T_X)-2\geq 0$.

Given a point $x\in X$, we denote by $H_x$ the normalization of the
subscheme of $H$ parametrizing rational curves passing through $x$.
By \cite[II.1.7, II.2.16]{kollar}, if $x\in X$ is a general point,
then $H_x$ is a smooth projective variety of dimension
$d=\deg(f^*T_X)-2$.  We remark that a rational curve that is smooth at
$x$ is parametrized by at most one element of $H_x$.

Let $H_1, \dots, H_k$ be minimal covering families of rational curves
on $X$.  For each $i$, let $\overline H_i$ denote the closure of $H_i$
in $\chow(X)$.  We define the following equivalence relation on $X$,
which we call $(H_1, \dots, H_k)$-equivalence.  Two points $x,y\in X$
are $(H_1, \dots, H_k)$-equivalent if they can be connected by a chain
of 1-cycles from $\overline H_1\cup \cdots \cup \overline H_k$.  By
\cite{campana} (see also \cite[IV.4.16]{kollar}), there exists a
proper surjective morphism $\pi^\circ:X^\circ \to Y^\circ$ from a
dense open subset of $X$ onto a normal variety whose fibers are $(H_1,
\dots, H_k)$-equivalence classes.  We call this map the \emph{$(H_1,
  \dots, H_k)$-rationally connected quotient of $X$}.  When $Y^\circ$
is a point we say that $X$ is $(H_1, \dots, H_k)$-rationally
connected.

\begin{rem}\label{limit_of_chains} 
  By \cite[IV.4.16]{kollar}, 
  there is a universal constant $c$, depending only on the dimension
  of $X$, with the following property.  If $H_1, \dots, H_k$ are
  minimal covering families of rational curves on $X$, and $x,y\in X$
  are general points on a general $(H_1, \dots, H_k)$-equivalence
  class, then $x$ and $y$ can be connected by a chain of  at most $c$ rational
  cycles from $\overline H_1\cup \cdots \cup \overline H_k$.
\end{rem}

The next two results are special features of the $(H_1, \dots,
H_k)$-rationally connected quotient of $X$ when the families $H_1,
\dots, H_k$ are unsplit.  The first one says that $\pi^\circ$ can be
extended in codimension $1$ to an equidimensional proper morphism with
integral fibers, but possibly allowing singular fibers.  The second
one describes the general fiber of the $H$-rationally connected
quotient of $X$ when $H$ is unsplit and $H_x$ is irreducible for
general $x\in X$.

\begin{lemma}\label{extending_in_codim_1}
  Let $X$ be a smooth complex projective variety and $H_1, \dots, H_k$
  unsplit covering families of rational curves on $X$.  Then there is
  an open subset $X^\circ$ of $X$, with $\codim_X(X\setminus
  X^\circ)\geq 2$, a smooth variety $Y^\circ$, and a proper surjective
  equidimensional morphism with irreducible and reduced fibers
  $\pi^\circ:X^\circ \to Y^\circ$ which is the $(H_1, \dots,
  H_k)$-rationally connected quotient of $X$.
\end{lemma}

\begin{proof}
  The fact that the $(H_1, \dots, H_k)$-rationally connected quotient
  of $X$ can be extended in codimension $1$ to an equidimensional
  proper morphism follows from the proof of \cite[Proposition
  1]{BCD_extremal_rays}.  This holds even in the more general context
  of quasi-unsplit covering families on $\q$-factorial varieties.
  In \cite[Proposition 1]{BCD_extremal_rays} this is established
  for a single quasi-unsplit family, but the same proof works for finitely many
  quasi-unsplit families. For convenience we review the construction of that extension.

  Let $\pi^\circ:X^\circ \to Y^\circ$ be the $(H_1, \dots,
  H_k)$-rationally connected quotient of $X$.  By shrinking $Y^\circ$
  if necessary, we may assume that $\pi^\circ$ is smooth.  Let $Y\to
  \chow(X)$ be the normalization of the closure of the image of
  $Y^\circ$ in $\chow(X)$, and let $\cU\subset Y\times X$ be the
  restriction of the universal family to $Y$. Denote by $p:\cU\to Y$
  and $q:\cU\to X$ the induced natural morphisms.  Notice that
  $q:\cU\to X$ is birational.

  Let $0\in Y$ and set $\cU_0=p^{-1}(0)$. Then $q(\cU_0)$ is contained
  in an $(H_1, \dots, H_k)$-equivalence class. This follows from
  taking limits of chains of rational curves from the families 
  $H_1, \dots, H_k$ (see Remark~\ref{limit_of_chains}), observing the
  assumption that the $H_i$'s are unsplit, and the fact that the image
  of a general fiber of $p$ in $X$ is an $(H_1, \dots,
  H_k)$-equivalence class.

  Let $E$ be the exceptional locus of $q$. Since $X$ is smooth, $E$
  has pure codimension $1$ in $\cU$. Set $S=q(E)\subset X$. This is a
  set of codimension at least $2$ in $X$. We shall show that $S$ is
  closed with respect to $(H_1, \dots, H_k)$-equivalence.  From that
  it will follow that the morphism $p|_{\cU\setminus E}:\cU\setminus
  E\to Y\setminus p(E)$ is proper and induces a proper equidimensional
  morphism $X\setminus S\to Y\setminus p(E)$ extending $\pi^\circ$.
  Let $L$ be an effective ample divisor on $Y$. Then there exists an
  effective $q$-exceptional divisor $F$ on $\cU$ and an effective
  divisor $D$ on $X$ such that $p^*L+F=q^*D$. First we claim that $\supp F =
  E$.  Indeed, let $C\subset E$ be any curve contracted by $q$. Then
  $C$ is not contracted by $p$ since $\cU\subset Y\times X$.  Hence
  $F\cdot C=q^*D\cdot C-p^*L\cdot C<0$, and so $C\subset \supp F$.
  This proves the claim.
  Notice that the general fiber of $p$ does not meet $E$.  Therefore,
  for any curve $C\subset \cU$ contained in a general fiber of $p$, we
  have $q^*D\cdot C=0$. This shows in particular that $D\cdot \ell=0$
  for any curve $\ell$ from any of the families $H_1, \dots, H_k$.  If
  $\tilde \ell\subset \cU$ is mapped onto $\ell$ by $q$, then $F\cdot
  \tilde \ell=q^*D\cdot \tilde \ell-p^*L\cdot \tilde \ell\leq 0$.
  Hence either $\tilde \ell$ is contained in $E=\supp F$ or it is
  disjoint from it.  Therefore, if $\ell$ is a curve from any of the
  families $H_1, \dots, H_k$, then either $\ell\subset S$ or $\ell\cap
  S=\emptyset$.  In other words, $S$ is closed with respect to $(H_1,
  \dots, H_k)$-equivalence.

  Replace $X^\circ$ with $X\setminus S$ and $Y^\circ$ with $Y\setminus
  p(E)$, obtaining a proper equidimensional morphism
  $\pi^\circ:X^\circ\to Y^\circ$ with $\codim(X\setminus X^\circ)\geq
  2$.  Since $Y$ is normal, we may also replace $Y^\circ$ with its
  smooth locus and we still have the condition $\codim(X\setminus
  X^\circ)\geq 2$.

  The locus $B$ of $Y^\circ$ over which $\pi^\circ$ has multiple
  fibers has codimension at least $2$ in $Y^\circ$.  To see this,
  compactify $Y^\circ$ to a smooth projective variety $\bar Y$ and
  take a resolution $\bar \pi :\bar X\to \bar Y$ of the
  indeterminacies of $X\map \bar Y$ with $\bar X$ smooth and
  projective.  Let $\bar C\subset \bar Y$ be a smooth projective curve
  obtained by intersecting $\dim \bar Y-1$ general very ample divisors
  on $\bar Y$. Let $\bar\pi_{\bar C}:\bar X_{\bar C}\to \bar C$ be the
  corresponding morphism.  Then $\bar X_{\bar C}$ is smooth projective
  and the general fiber of $\bar \pi$ is rationally connected.  Hence
  $\bar\pi_{\bar C}$ has a section by \cite{GHS}, and thus it cannot
  contain multiple fibers. Now, replace $Y^\circ$ with
  $Y^\circ\setminus B$ to obtain an equidimensional proper morphism
  with no multiple fibers.

  Let $F$ be a general fiber of $\pi^\circ$. For each $i$, denote by
  $H_i^j$, $1\leq j\leq n_i$, 
  the unsplit covering families of
  rational curves on $F$ whose general members correspond to rational 
  curves on $X$ from the family $H_i$. Let
  $[H_i^j]$ denote the class of a member of $H_i^j$ in $N_1(F)$ and
  $\cH:=\{[H_i^j]\mid i=1,\dots,k, j=1,\dots,n_i\}$.  Then by
  \cite[IV.3.13.3]{kollar}, $N_1(F)$ is generated by $\cH$.

  Finally we shall show that the locus $B'$ of $Y^\circ$ over which
  the fibers of $\pi^\circ$ are not integral has codimension at least
  $2$ in $Y^\circ$.  Let $C\subset Y^\circ$ be a smooth curve obtained
  by intersecting $\dim Y^\circ-1$ general very ample divisors on
  $Y^\circ$.  Let $\pi_{C}:X_{C}\to C$ be the corresponding morphism.
  Then $X_{C}$ is smooth.  We denote the image of the classes
  $[H_i^j]$'s in $N_1(X_C)$ and their collection $\cH$ by the same
  symbols.  By taking limits of chains of rational curves from 
  the families $H_1, \dots, H_k$ and
  applying \cite[IV.3.13.3]{kollar} (see Remark~\ref{limit_of_chains}), 
  we see that any curve contained
  in any fiber of $\pi_C$ is numerically proportional in $N_1(X_C)$ to
  a linear combination of the $[H_i^j]$'s.  Hence $N_1(X_C/C)$ is
  generated by $\cH$.  Therefore, all fibers of $\pi_C$ are
  irreducible.  Indeed, if $F_0'$ is an irreducible component of a
  reducible fiber $F_0$, then $F_0'$ is a Cartier divisor on $X_C$,
  and $F_0'\cdot [H_i^j]=0$ for every $H_i^j$.  On the other hand,
  there is a curve $\ell \subset F_0$ such that $F_0'\cdot \ell >0$,
  contradicting the fact that $N_1(X_C/C)$ is generated by $\cH$.
  Since there are no multiple fibers, the fibers are also reduced.
  Finally, we replace $Y^\circ$ with $Y^\circ\setminus B'$ and obtain
  a morphism with the required properties.
\end{proof}

\begin{prop}\label{fibers_have_rho=1}
  Let $X$ be a smooth complex projective variety and $H$ an unsplit
  covering family of rational curves on $X$. Assume that $H_x$ is
  irreducible for general $x\in X$.  Let $\pi^\circ:X^\circ \to Y^\circ$
  be the $H$-rationally connected quotient of $X$.  Then the general fiber of
  $\pi^\circ$ is a Fano manifold with Picard number $1$.
\end{prop}

\begin{proof}
  Let $X_t$ be a general fiber of $\pi^\circ$, and suppose
  $\rho(X_t)\neq 1$.  Denote by $[H]$ the class of the members of $H$
  in $N_1(X)$.  By \cite[IV.3.13.3]{kollar}, every proper curve on
  $X_t$ is numerically proportional to $[H]$ in $N_1(X)$.  There
  exists an irreducible component $H_t$ of $H_{X_t}=\{[C]\in H\ |\
  C\subset X_t\}$ which is an unsplit covering family of rational
  curves on $X_t$.  Since $H_x$ is irreducible for general $x\in X$,
  such a component $H_t$ is unique.  Since $\rho(X_t)\neq 1$, $X_t$ is
  not $H_t$-rationally connected by \cite[IV.3.13.3]{kollar}.  Let
  $\sigma_t:X_t^\circ \to Z_t^\circ$ be the (nontrivial)
  $H_t$-rationally connected quotient of $X_t$.  Notice that for every
  $z\in Z_t^\circ$ there is a curve $C_z\subset X_t$ numerically
  proportional to $[H]$ in $N_1(X)$, meeting the fiber of $\sigma_t$
  over $z$, but not contained in it.  Since $H_t$ is unique, there is
  a dense open subset $X'$ of $X$ and a fibration $\sigma:X' \to Z'$
  whose fibers are fibers of $\sigma_t$ for some $t\in Y^\circ$.
  Moreover, there is a curve $C\subset X$ numerically proportional to
  $[H]$ in $N_1(X)$, meeting $X'$, and not contracted by $\sigma$.
  But this is impossible.  Indeed, let $L'$ be an effective divisor on
  $Z'$ meeting but not containing the image of $C$ by $\sigma$.  Let
  $L$ be the closure of $\sigma^{-1}(L')$ in $X$.  Then $L\cdot C>0$
  while $L\cdot \ell=0$ for any curve $\ell$ parametrized by $H$ lying
  on a fiber of $\sigma$.
\end{proof}

\begin{rem}
  The statement of Proposition~\ref{fibers_have_rho=1} does not hold
  in general if we do not assume that $H_x$ is irreducible for general
  $x\in X$.  Indeed, one may take $\pi^\circ:X^\circ\to Y^\circ$ to be
  a suitable family of quadric surfaces in $\p^3$ and $H$ to be the
  family of lines on the fibers of $\pi^\circ$.
\end{rem}

\begin{defn}
  Let $X$ be a smooth complex projective variety, and $H$ a minimal
  covering family of rational curves on $X$.  Let $x\in X$ be a
  general point. Define the tangent map $ \tau_x: \ H_x \map
  \bP(T_xX^*) $ by sending a curve that is smooth at $x$ to its
  tangent direction at $x$.  Define $\cC_x$ to be the image of
  $\tau_x$ in $\bP(T_xX^*)$.  This is called the \emph{variety of
    minimal rational tangents} at $x$ associated to the minimal family
  $H$.
\end{defn}

The map $\tau_x: \ H_x \to \cC_x$ is in fact the normalization
morphism by \cite{kebekus} and \cite{hwang_mok_birationality}.  If
$\tau_x$ is an immersion at every point of $H_x$, then all curves
parametrized by $H_x$ are smooth at $x$ by \cite[V.3.6]{kollar} and
\cite[Proposition 2.7]{artigo_tese}, and, as a consequence, 
there is a one-to-one corresponcence between points of $H_x$ 
and the associated curves on $X$.
The variety $\cC_x$ comes with a natural projective embedding into
$\bP(T_xX^*)$.  This embedding encodes important geometric properties
of $X$. The following result was proved in \cite{artigo_tese} and
gives a structure theorem for varieties whose variety of minimal
rational tangents is linear.

\begin{thm}[\cite{artigo_tese}]\label{pn_bundle_in_codim_1} 
  Let $X$ be a smooth complex projective variety, $H$ a minimal
  covering family of rational curves on $X$, and $\cC_x\subset
  \bP(T_xX^*)$ the corresponding variety of minimal rational tangents
  at $x\in X$. Suppose that for a general $x\in X$, $\cC_x$ is a
  $d$-dimensional linear subspace of $\bP(T_xX^*)$.
  
  Then there exists an open subset $X^\circ\subset X$ and a
  $\bP^{d+1}$-bundle $\varphi^\circ:X^\circ\to T^\circ$ over a smooth
  base with the property that every rational curve parametrized by $H$
  and meeting $X^\circ$ is a line on a fiber of $\varphi^\circ$.  In
  particular, $\varphi^\circ:X^\circ\to T^\circ$ is the $H$-rationally
  connected quotient of $X$.
  If $H$ is unsplit, then we may take $X^\circ$ such that
  $\codim(X\setminus X^\circ)\geq 2$.
\end{thm}

\begin{prop}\label{ample_subsheaf_of_f*TX}
  Let $X$ be a smooth complex projective variety, $H$ a minimal
  covering family of rational curves on $X$, and $\pi^\circ:X^\circ
  \to Y^\circ$ the $H$-rationally connected quotient of $X$.  Suppose
  that the tangent bundle $T_X$ contains a subsheaf $\sD$ such that
  $f^*\sD$ is an ample vector bundle for a general member $[f]\in H$.
  Then $\pi^\circ$ is a projective space bundle and the inclusion
  $\sD|_{X^\circ}\into T_{X^\circ}$ factors through an inclusion
  $\sD|_{X^\circ}\into T_{X^\circ/Y^\circ}$.
\end{prop}

\begin{proof}
  Let $\cC_x\subset \bP(T_xX^*)$ be the variety of minimal rational
  tangents associated to $H$ at a general point $x\in X$.  By
  \cite[Proposition 4.1]{artigo_tese}, $\cC_x$ is a union of linear
  subspaces of $\bP(T_xX^*)$ containing $\bP(\sD^*\otimes \kappa(x))$.
  In \cite[Proposition 4.1]{artigo_tese} $\sD$ is assumed to be ample,
  but the proof only uses the fact that $f^*\sD$ is a subsheaf of
  $f^*T_X^+$ for general $[f]\in H$.

 Lemma~\ref{Cx_must_be_irred} below implies that $\cC_x$ is irreducible, 
 and thus a linear subspace of $\bP(T_xX^*)$.

  Now we apply Theorem~\ref{pn_bundle_in_codim_1} to conclude that
  $\pi^\circ$ is a projective space bundle.  Moreover, for a general
  point $x\in X^\circ$, the stalk $\sD_x$ is contained in
  $(T_{X^\circ/Y^\circ})_x$.  Since the cokernel of
  $T_{X^\circ/Y^\circ}\into T_{X^\circ}$ is torsion free, we conclude
  that there is an inclusion $\sD|_{X^\circ}\into T_{X^\circ/Y^\circ}$
  factoring $\sD|_{X^\circ}\into T_{X^\circ}$.
\end{proof}

The following lemma is Proposition 2.2 of \cite{hwang_b2=b4=1}.  In
\cite[Proposition 2.2]{hwang_b2=b4=1} $X$ is assumed to have Picard
number one, but this assumption is not used in the proof.

\begin{lemma}\label{Cx_must_be_irred}
  Let $X$ be a smooth complex projective variety, $H$ a minimal
  covering family of rational curves on $X$, and $\cC_x\subset
  \bP(T_xX^*)$ the corresponding variety of minimal rational tangents
  at $x\in X$. Suppose that for a general $x\in X$, $\cC_x$ is a
  union of linear subspaces of $\bP(T_xX^*)$.

  Then the intersection of any two irreducible components of $\cC_x$ is empty. 

\end{lemma}



\section{The relative anticanonical bundle of a fibration}
\label{section_relative_-K}

In this section we prove that the relative anticanonical bundle of a
generically smooth surjective morphism from a normal projective
$\q$-Gorenstein variety onto a smooth curve cannot be ample.  In fact,
we prove the following more general result.  Note that a similar
theorem was proved in \cite[Theorem 2]{miyaoka_relative_deformations}.

\begin{thm}\label{-KX/Y_not_ample}
  Let $X$ be a normal projective variety, $f:X\to C$ a surjective morphism onto a smooth
  curve, and $\Delta\subseteq X$ a Weil
  divisor such that $(X,\Delta)$ is log canonical over the generic
  point of $C$. Then $-(K_{X/C}+\Delta)$ is not ample.
\end{thm}

\begin{proof}
  Let $X\overset g\to \tilde C\overset\sigma\to C$ be the Stein
  factorization of $f$.  Then $K_{\tilde C}=\sigma^*K_C+R_\sigma$
  where $R_\sigma$ is the ramification divisor of $\sigma$ and so
  $-(K_{X/\tilde C}+\Delta)=-(K_{X/C}+\Delta)+{g}^*R_\sigma$.  Notice
  that $R_\sigma$ is effective and hence if $-(K_{X/C}+\Delta)$ is
  ample, then so is $-(K_{X/\tilde C}+\Delta)$.
  
  Thus, in order to prove the statement, we may assume that $f$ has
  connected fibers. Let us now assume to the contrary that
  $-(K_{X/C}+\Delta)$ is ample.  Let $\pi: \tilde X\to X$ be a log
  resolution of singularities of $(X,\Delta)$, $A$ an ample divisor on
  $C$, and $m\gg 0$ such that $D=-m(K_{X/C}+\Delta)-f^*A$ is very
  ample.  Then
  $$
  K_{\tilde X}+\pi^{-1}_* \Delta\sim_{\q} \pi^*(K_X +\Delta)+E_+-E_-,
  $$
  where $E_+$ and $E_-$ are effective $\pi$-exceptional divisors with
  no common components and such that the support of
  $\pi^{-1}_*\Delta+E_++E_-$ is an snc divisor.  By the log canonical
  assumption, $E_-$ can be decomposed as $E_-=E+F$ where $\rup E$ is
  reduced and $E_-$ agrees with $E$ over the generic point of $C$.
  Set $\tilde f=f\circ \pi$ and let $\tilde D\in |\pi^*D|$ be a general
  member. Setting $\tilde\Delta=\pi^{-1}_*\Delta+\frac 1m\tilde D+E$,
  we obtain that $(\tilde X,\tilde \Delta)$ is log canonical and that
  \begin{equation}
    \label{eq:1}
    K_{\tilde X}+\tilde \Delta + F \sim_{\q} \tilde f^*K_C +E_+
    -\frac{1}{m}\tilde f^*A. 
  \end{equation}
 
  Furthermore, since $E_+$ is $\pi$-exceptional, $\pi_*\sO_{\tilde
    X}(lE_+)$ is an ideal sheaf in $\sO_X$ for any $l\in\bZ$ (see for
  instance \cite[Lemma 7.11]{debarre}).  Then for any $l\in \bN$
  sufficiently divisible,
  \begin{multline*}
    \tilde f_*\sO_{\tilde X}(lm(K_{\tilde X/C}+\tilde \Delta))
    \overset\iota\into \tilde f_*\sO_{\tilde X}(lm(K_{\tilde
      X/C}+\tilde
    \Delta+F))\simeq  \\
    \simeq \tilde f_*\sO_{\tilde X}(l(mE_+-\tilde f^*A)) \simeq \tilde
    f_*\sO_{\tilde X}(lmE_+)\otimes \sO_C(-lA) \subseteq \sO_C(-lA).
  \end{multline*}
  Finally, observe that 
  \begin{itemize}
  \item $\tilde f_*\sO_{\tilde X}(lm(K_{\tilde X/C}+\tilde \Delta+F))$
    is nonzero by \eqref{eq:1} and because $E_+$ is effective,
  \item $\tilde f_*\sO_{\tilde X}(lm(K_{\tilde X/C}+\tilde\Delta))$ is
    semi-positive by \cite[Thm.~4.13]{campana_orbifold}, and
  \item $\iota$ is an isomorphism over a nonempty open subset of $C$.
  \end{itemize}
  Therefore, $\tilde f_*\sO_{\tilde X}(lm(K_{\tilde
    X/C}+\tilde\Delta))$ is a non-zero semi-positive sheaf contained
  in $\sO_C(-lA)$, but that contradicts the fact that $A$ is ample.
\end{proof}



\section{Lifting $p$-derivations to the normalization}
\label{section_lifting_p_derivations}

In this section we show that $p$-derivations (see
definition~\ref{def_p_derivation_on_X} below) can be lifted to the
normalization.  This is a generalization of Seidenberg's theorem in
\cite{seidenberg}.  The proofs in this section follow closely the
proof of Theorem 2.1.1 in \cite{kallstrom} and we also use the
following result from \cite{kallstrom}.

\begin{lemma}[{\cite[Lemma 2.1.2]{kallstrom}}]\label{kallstrom} 
  Let $(A,\m,k)$ be a local Noetherian domain and $\partial$ a
  derivation of $A$. Let $\nu$ be a discrete valuation on the fraction
  field $K(A)$ with center in $A$.  Then there exists a $c\in\z$ such
  that $\nu\left(\frac{\partial(x)}{x}\right)\geq c$ for any $x\in
  K(A)\setminus\{0\}$.
\end{lemma}

\begin{defn}\label{def_p_derivation_on_A}
  Let $R$ be a ring, $A$ an $R$-algebra and $M$ an $A$-module.  Denote
  by $\Omega_{A/R}$ the module of relative differentials of $A$ over
  $R$.  Given a positive integer $p$, we denote by $\Omega_{A/R}^p$
  the $p$-th wedge power of $\Omega_{A/R}$.  A \emph{$p$-derivation of
    $A$ over $R$ with values in $M$} is an $A$-linear map
  $\partial:\Omega_{A/R}^p\to M$. Such a map $\partial$ induces a skew
  symmetric map $K(A)^{\oplus p}\to M\otimes_A K(A)$, where $K(A)$
  denotes the fraction field of $A$.  We use the same symbol
  $\partial$ to denote this induced map.  When $M=A$ and $R$ is clear
  from the context, we call $\partial$ simply a \emph{$p$-derivation
    of $A$}.
\end{defn}

\begin{lemma}\label{continuity}
  Let $(A,\m,k)$ be a local Noetherian domain, $p$ a positive integer,
  and $\partial$ a $p$-derivation of $A$.  Let $\nu$ be a discrete
  valuation on the fraction field $K(A)$ with center in $A$.  Then
  there exists $c\in\z$ such that $\nu\left(\frac{\partial(x_1,\dots,
      x_p)}{x_1 \cdots x_p}\right)\geq c$ for any $x_1,\dots, x_p\in
  K(A)\setminus\{0\}$.
\end{lemma}

\begin{proof}
  We use induction on $p$.  If $p=1$, this is Lemma~\ref{kallstrom}.
  Suppose now that $p \geq 2$ and let $(A,\m,k)$ be a local Noetherian
  domain, $\partial$ a $p$-derivation of $A$, and $\nu$ a discrete
  valuation on the fraction field $K(A)$ with center in $A$.  Let
  $m_{1},\dots,m_{r}$ be generators of the maximal ideal $\m$.
  
  Using the formula
  $$
  \frac{\partial(x_{1,1}x_{1,2},\ldots,x_{p,1}x_{p,2})}
  {x_{1,1}x_{1,2}\cdots x_{p,1}x_{p,2}} =\sum
  \frac{\partial(x_{1,i_{1}},\ldots,x_{p,i_{p}})}{x_{1,i_{1}}\cdots
    x_{p,i_{p}}},
  $$
  we get
  $$
  \nu\left(\frac{\partial(x_{1,1}x_{1,2},\ldots,x_{p,1}x_{p,2})}
    {x_{1,1}x_{1,2}\cdots x_{p,1}x_{p,2}}\right)
  \ge\min\left\{\nu\left(
      \frac{\partial(x_{1,i_{1}},\ldots,x_{p,i_{p}})}{x_{1,i_{1}}\cdots
        x_{p,i_{p}}}\right)\right\}
  $$
  for $x_{1,1},x_{1,2},\ldots,x_{p,1},x_{p,2}\in A\setminus\{0\}$.
  Further observe that
  $$
  \frac{\partial(x_{1}^{-1},x_{2},\ldots,x_{p})} {x_{1}^{-1}x_{2}\cdots
    x_{p}}=-\frac{\partial(x_{1},\ldots,x_{p})}{x_{1}\cdots x_{p}}.
  $$
  Also, if $a\in A$, then $a$ may be written as a sum of products
  $m_{i_{1}}\cdots m_{i_{k}}u$ with $u\in A\setminus\m$.  Therefore we
  only have to check that the required inequality holds for
  $x_{1},\ldots,x_{p}\in\{m_{1},\ldots,m_{r}\}\cup(A\setminus\m)$.

  If $x_{1},\ldots,x_{p}\in A\setminus\m$ then
  $$
  \nu\left(\frac{\partial(x_{1},\ldots,x_{p})}{x_{1}\cdots
      x_{p}}\right)= \nu(\partial(x_{1},\ldots,x_{p}))\ge 0.
  $$

  Suppose now that at least one of the $x_i$'s is in $\m$.  For
  simplicity we assume that $x_{1},\ldots,x_{l}\in A\setminus\m$ and
  $x_{l+1},\ldots,x_{p}\in\{m_{1},\ldots,m_{r}\}$, $0\leq l< p$.  We
  may view $\partial(\cdot,\cdots,\cdot,x_{l+1},\ldots,x_{p})$ as an
  $l$-derivation of $A$.  The result then follows by induction.
\end{proof}

\begin{defn}\label{def_p_derivation_on_X}
  Let $S$ be a scheme, $X$ a scheme over $S$, and $\sL$ a line bundle
  on $X$.  Denote by $\Omega_{X/S}$ the sheaf of relative
  differentials of $X$ over $S$, and by $\Omega_{X/S}^p$ its $p$-th
  wedge power for $p\in\bN$.  A \emph{$p$-derivation of $X$ over $S$
    with values in $\sL$} is a morphism of sheaves
  $\partial:\Omega_{X/S}^p\to \sL$.  When $S$ is the spectrum of a
  field and $\sL$ is clear from the context, we drop $S$ and $\sL$
  from the notation and call $\partial$ simply a \emph{$p$-derivation
    on $X$}.
\end{defn}

\begin{prop}\label{extensionpderivation}
  Let $X$ be a Noetherian integral scheme over a field $k$ of
  characteristic zero and $\eta:\widetilde X\to X$ its normalization.  Let
  $\sL$ be a line bundle on $X$, $p$ a positive integer, and
  $\partial$ a $p$-derivation with values in $\sL$.  Then $\partial$
  extends to a unique $p$-derivation $\bar\partial$ on $\widetilde X$ with
  values in $\eta^{*}\sL$.
\end{prop}

\begin{proof}
  The uniqueness of $\bar\partial$ is clear since $\sL$ is torsion
  free and $\eta$ is birational.  The existence of the lifting can be
  established locally.  So we may assume that $X$ is the spectrum of
  an integral $k$-algebra $A$, $\sL$ is trivial, and $\partial$ is a
  $p$-derivation of $A$.  Let $A'$ denote the integral closure of $A$
  in its fraction field $K(A)$.  There exists a unique extension of
  $\partial$ to a $p$-derivation of $K(A)$, which we also denote by
  $\partial$.  We must prove that $\partial(A',\ldots,A')\subset A'$.

  First we reduce the problem to the case when $A$ is a
  $1$-dimensional local ring and $A'$ is a DVR.  Since $A'$ is
  integrally closed in $K(A)$, $A'$ is the intersection of its
  localizations at primes of height one \cite[2. Theorem
  38]{matsumura}.  Let $\mathfrak{p}'$ be a prime of height one of
  $A'$, and set $\mathfrak{p}=\mathfrak{p}'\cap A$.  Notice that
  $\partial(A_{\mathfrak{p}},\ldots,A_{\mathfrak{p}})\subset
  A_{\mathfrak{p}}$, and the result follows if we prove that
  $\partial(A'_{\mathfrak{p'}},\ldots,A'_{\mathfrak{p'}})\subset
  A'_{\mathfrak{p'}}$.  Hence we may assume that $A$ is a
  $1$-dimensional local ring and $A'$ is a DVR.  Denote by $\m$ and
  $\m'$ the maximal ideals of $A$ and $A'$ respectively.

  Next we further reduce the problem to the case when $A$ and $A'$ are
  complete local rings.  Let $\bar R$ be the completion of $A'$ with
  respect to the $\m'$-adic topology.  Let $\bar A$ be the completion
  of $A$ with respect to the $\m$-adic topology. Since $\bar A$ is
  $1$-dimensional, there is an inclusion of local rings $\bar A\subset
  \bar R$.  Let $\nu$ be a discrete valuation of $K(A')$ whose
  valuation ring is $A'$.  By Lemma~\ref{continuity}, $\partial$ is a
  continuous $p$-derivation of $R$ with values in $K(A')$.  Hence it
  has a unique extension to a continuous $p$-derivation $\bar
  \partial$ of $K(\bar{R})$.  Notice that the condition
  $\partial(A,\ldots,A)\subset A$ implies that
  $\partial(\bar{A},\ldots,\bar{A})\subset\bar{A}$ by the Artin-Rees
  Lemma.  Since $K(A)\cap\bar R=A'$, the result then follows if we
  prove that $\partial(\bar{R},\ldots,\bar{R})\subset\bar{R}$.
  Therefore we may assume that $A$ and $A'$ are complete
  $1$-dimensional local rings.

  Now we use induction on $p$. If $p=1$, this is Seidenberg's theorem
  \cite{seidenberg}, so we may assume that $p\geq 2$. Let $k_{A}$ be a
  coefficient field in $A$, and $k_{A'}$ a coefficient field in $A'$
  containing $k_{A}$ \cite[Theorem 7.8]{eisenbud}.  The extension
  $k_{A'}|k_{A}$ is finite.  Let $t\in\m'$ be a uniformizing
  parameter.  It suffices to show that
  $\partial(x_{1},\ldots,x_{p})\in A'$ for $x_{1},\ldots,x_{p}\in
  k_{A'}\cup\{t\}$.  Since $\partial$ is skew symmetric and $p\ge 2$,
  we have $\partial(t,\ldots,t)=0$.  So we may assume that $x_{1}\in
  k_{A'}$. Since $k_{A'}|k_{A}$ is finite and separable, there exists
  $P(X)=\sum a_{i} X^{i}\in k_{A}[X]$ such that $P(x_{1})=0$ and
  $P'(x_{1})\neq 0$.  Thus
  $$
  0=\partial(P(x_{1}),x_{2},\ldots,x_{p})
  =P'(x_{1})\partial(x_{1},\ldots,x_{p}) + \sum\partial(a_{i},x_{2},
  \ldots,x_{p})x_{1}^{i}.
  $$
  Finally, $\partial(a_{i},\_,\ldots,\_)$ may be viewed as a $p-1$
  derivation of $A$ and so $\partial(x_{1},\ldots,x_{p})\in A'$ by the
  induction hypothesis.
\end{proof}



\section{Sections of $\wedge^{p}T_{X}\otimes \sM$}
\label{section_fibration}

The following lemma will be used several times in this section.

\begin{lemma}\label{2_cases}
  Let $Y$ be a smooth variety, $\pi : X\to Y$ a smooth morphism, $\sM$
  a line bundle on $X$, and $p\geq 2$ an integer.  Suppose that for a
  general fiber, $F$, of $\pi$, $H^0(F,\wedge^{i}T_{F}\otimes \sM|_F)=
  0$ for $0\leq i\leq p-2$.  Then there exists an exact sequence:
  \begin{multline*}
    0 \to H^0(X,\wedge^{p}T_{X/Y}\otimes \sM) \to
    H^0(X,\wedge^{p}T_{X}\otimes \sM) \to 
    H^0(X,\wedge^{p-1}T_{X/Y}\otimes\pi^{*}{T_{Y}}\otimes \sM).
  \end{multline*}
\end{lemma}

\begin{proof}
  The short exact sequence
  $$
  0\to T_{X/Y}\to T_{X}\to \pi^{*}T_{Y}\to 0
  $$
  yields a filtration $\wedge^{p}T_{X}\otimes \sM=\sF_{0} \supseteq
  \sF_{1}\supseteq \sF_{2}\supseteq \cdots \supseteq \sF_p \supseteq
  \sF_{p+1}=0$ such that
  $$
  \sF_{i}/\sF_{i+1}\simeq \wedge^{i}
  T_{X/Y}\otimes\pi^{*}\wedge^{p-i}{T_{Y}}\otimes \sM
  $$
  for each $i$. In particular, one has the short exact sequence,
  \begin{equation}
    \label{eq:2}
    0\to \wedge^{p}T_{X/Y}\otimes \sM \to \sF_{p-1} \to
    \wedge^{p-1}T_{X/Y}\otimes\pi^{*}{T_{Y}}\otimes \sM \to 0. 
  \end{equation}
  The assumption that $H^0(F,\wedge^{i}T_{F}\otimes \sM|_F)= 0$ for
  $0\leq i\leq p-2$ for a general fiber of $\pi$ implies that
  $H^0(X,\sF_{i}/\sF_{i+1})=0$ for $0\leq i\leq p-2$, thus $H^0(X,
  \wedge^{p}T_{X}\otimes \sM)= H^0(X,\sF_{0})=\dots=H^0(X,\sF_{p-1})$
  and the result follows from (\ref{eq:2}).
\end{proof}

The condition that $H^0(F,\wedge^{i}T_{F}\otimes \sM|_F)= 0$ for $0\leq
i\leq p-2$ and $F$ a general fiber of $\pi$ is easily verified when
$\pi$ is a projective space bundle and $\sM|_F$ is sufficiently
negative. In this case we get the following.

\begin{lemma}\label{lemma:pbundle}
  Let $Y$ be a smooth projective variety of dimension $\ge 1$, $\sE$
  an ample vector bundle of rank $r+1\geq 2$ and $\sN$ a nef line
  bundle on $Y$.  Consider the projective bundle $\pi:X=\bP(\sE)\to Y$
  with tautological line bundle $\sO_{\bP(\sE)}(1)$.  Let $p,q\in\bN$
  and assume that $p\geq 2$. Then
  \begin{equation}
    \label{eq:5}
    H^0(X,\wedge^{p}T_{X/Y}\otimes \sO_{\bP(\sE)}(-p-q)\otimes
    \pi^*\sN^{-1})= 0.
  \end{equation}
\end{lemma}

\begin{proof} 
  First observe, that if $p>r$ then the statement is trivially true,
  so we will assume that $p\leq r$. Let $i\in\bN$, $i<p$.  After
  twisting by $\sO_{\bP(\sE)}(-p-q)\otimes \pi^*\sN^{-1}$, the short
  exact sequence
  $$
  0\to \wedge^{p-i-1}T_{X/Y} \to
  \wedge^{p-i}\left(\pi^{*}\sE^{\dual}(1)\right)\to
  \wedge^{p-i}T_{X/Y}\to 0
  $$
  yields the exact sequence
    \begin{multline}
      \label{eq:3}
      \dots \to H^{i}(X,\wedge^{p-i}(\pi^{*}\sE^{\dual})(-i-q)\otimes
      \pi^*\sN^{-1}) 
      \to H^{i}(X,\wedge^{p-i}T_{X/Y}(-p-q)\otimes \pi^*\sN^{-1})\to
      \\ 
      \to H^{i+1}(X,\wedge^{p-i-1}T_{X/Y}(-p-q)\otimes \pi^*\sN^{-1})
      \to \dots\hfil
    \end{multline}
  Since $i<p\leq r$ and $R^j\pi_*\sO_{\bP(\sE)}(l)=0$ for $0<j<r$ and
  for any $l\in\bZ$, the Leray spectral sequence implies that
  $$
  H^{i}(X,\wedge^{p-i}(\pi^{*}\sE^{\dual})(-i-q)\otimes
  \pi^*\sN^{-1})= H^{i}(Y,\wedge^{p-i}\sE^{\dual}\otimes
  \sN^{-1}\otimes \pi_{*}\sO_{\bP(\sE)}(-i-q)).
  $$
  The sheaf $\pi_{*}\sO_{\bP(\sE)}(-i-q)$ is zero unless $i=q=0$, in which
  case it is isomorphic to $\sO_{Y}$.  Furthemore,
  $H^{0}(Y,\wedge^{p}\sE^{\dual}\otimes \sN^{-1})=0$ since $\sE$ is
  ample and $\sN$ is nef, and hence
  $$
  H^{i}(X,\wedge^{p-i}(\pi^{*}\sE^{\dual})(-i-q)\otimes
  \pi^*\sN^{-1})=0
  $$
  for $0\leq i\leq p-1$. Therefore, by (\ref{eq:3}), one has a series
  of injections,
  \begin{multline*}
    H^{0}(X,\wedge^{p}T_{X/Y}(-p-q)\otimes \pi^*\sN^{-1}) \into 
    H^{1}(X,\wedge^{p-1}T_{X/Y}(-p-q)\otimes \pi^*\sN^{-1}) \into
    \dots \\ \dots \into 
    H^{i}(X,\wedge^{p-i}T_{X/Y}(-p-q)\otimes \pi^*\sN^{-1}) 
    \into \dots \into 
    H^{p}(X,\sO_{\bP(\sE)}(-p-q)\otimes \pi^*\sN^{-1}).
  \end{multline*}
  
  By the Kodaira vanishing theorem
  $H^{p}(X,\sO_{\bP(\sE)}(-p-q)\otimes \pi^*\sN^{-1})= 0$, and the
  statement follows.
\end{proof}

\begin{cor}\label{pbundle}
  Let $Y$ be a smooth projective variety of dimension $\ge 1$ and
  $\sE$ an ample vector bundle of rank $r+1\ge 2$ on $Y$.  Consider
  the projective bundle $\pi:X=\bP(\sE)\to Y$ with tautological line
  bundle $\sO_{\bP(\sE)}(1)$.  Suppose that $H^0(X,\wedge^{p}T_{X}
  \otimes \sO_{\bP(\sE)}(-p-q)\otimes\pi^*\sN^{-1})\neq 0$ for some
  integers $p\ge 2$, $q\geq 0$, and some nef line bundle $\sN$ on $Y$.
  Then $Y\simeq \bP^1$, $\sE\simeq \sO_{\bP^1}(1)\oplus \sO_{\bP^1}(1)$,
  $p=2$, $q=0$, and $\sN\simeq\sO_{\bP^1}$.
\end{cor}

\begin{proof}
  Let $F\simeq \bP^{r}$ denote a general fiber of $\pi$ and set
  $\sM=\sO_{\bP(\sE)}(-p-q)\otimes \pi^*\sN^{-1}$. Then by Bott's
  formula $H^0(F,\wedge^{i}T_{F}\otimes \sM|_F)= 0$ for every $0\leq
  i\leq p-2$.  Then Lemma~\ref{2_cases} and Lemma~\ref{lemma:pbundle}
  imply that $H^0(X,\wedge^{p-1}T_{X/Y}\otimes\pi^{*} ({T_{Y}}\otimes
  \sN^{-1})\otimes \sO_{\bP(\sE)}(-p-q))\neq 0$. By Bott's formula
  again $H^0(F,\wedge^{p-1}T_{F}(-p-q))\neq 0$ implies that $q=0$ and
  $r=p-1$.  Therefore we have
  \begin{multline*}
    0\neq \coh 0.X.{\wedge^{r}T_{X/Y}\otimes\pi^{*}({T_{Y}}\otimes
      \sN^{-1})\otimes \sO_{\bP(\sE)}(-r-1)}.=\\ \coh 0. X.
    {\pi^*\left({T_{Y}}\otimes \det\sE^\dual\otimes
        \sN^{-1}\right)}.  \simeq \coh 0. Y.
    {\pi_*\pi^*\left({T_{Y}}\otimes \det \sE^{\dual}\otimes
        \sN^{-1}\right)}. \simeq \\ \coh 0. Y.  {{T_{Y}}\otimes
      \big(\underbrace{\det \sE\otimes
        \sN}_{\text{ample}}\big)^{-1}}..
  \end{multline*}
  Now Wahl's theorem \cite{wahl} yields that $Y\simeq \bP^m$ for some
  $m>0$. Then we immediately obtain that $\deg(\det \sE\otimes
  \sN)\leq 2$. Since $\sE$ is ample on a projective space,
  $$
  2\leq r+1=\rk\sE\leq \deg\sE \leq \deg(\det \sE\otimes \sN)
  -\deg\sN\leq 2- \deg\sN\leq 2.
  $$
  Therefore all of these inequalities must be equalities and we have
  that $r+1=p=2$, $q=0$ and $\sN\simeq \sO_Y$. Furthermore, this
  implies that then $\sO_{\bP^m}(2)\simeq \det\sE\into T_{\bP^m}$ and
  hence $m=1$.
\end{proof}


\begin{prop}\label{generalvanishing}
  Let $X$ be a smooth projective variety, $H\subset\rat(X)$ a minimal
  covering family of rational curves on $X$, $\sL$ an ample line
  bundle on $X$, and $\sM$ a nef line bundle on $X$ such that
  $c_1(\sM)\cdot C>0$ for every $[C]\in H$.  Suppose that
  $H^0(X,\wedge^{p} T_{X}\otimes \sL^{-p}\otimes \sM^{-1})\neq 0$ for
  some integer $p\ge 1$.  Then $(X,\sL,\sM)\simeq
  (\bP^{p},\sO_{\bP^{p}}(1),\sO_{\bP^{p}}(1))$.
\end{prop}

\begin{proof}
  Let $[f]\in H$ be a general member and write $f^{*}T_{X}\simeq
  \sO_{\bP^1}(2)\oplus\sO_{\bP^1}(1)^{\oplus
    d}\oplus\sO_{\bP^1}^{\oplus n-d-1}$.  The condition that both $f^*\sL$ and  $f^*\sM$
  are ample and that $H^0(X,\wedge^{p}T_{X}\otimes \sL^{-p}\otimes
  \sM^{-1})\neq 0$ implies that $f^*\sL\simeq \sO_{\bP^1}(1) \simeq
  f^*\sM$, and thus $H$ is unsplit.  A non-zero section $s\in
  H^0(X,\wedge^{p} T_{X}\otimes \sL^{-p}\otimes \sM^{-1})$ and the
  contraction 
  $$
  \sC_\theta:\wedge^pT_X\otimes \sL^{-p}\otimes \sM^{-1}\to
  \wedge^{p-1}T_X\otimes \sL^{-p}\otimes \sM^{-1}
  $$
  induced by a differential form $\theta\in\Omega_X$, gives rise to a
  non-zero map
  \begin{align*}
    \Omega_X &\to \wedge^{p-1}T_X\otimes \sL^{-p}\otimes \sM^{-1}\\
    \theta &\mapsto \sC_\theta(s),
  \end{align*}
  the dual of which is the non-zero map
  \begin{equation}
    \label{eq:4}
      \varphi:\Omega_{X}^{p-1}\otimes \sL^{p}\otimes \sM\to T_{X}.
  \end{equation}
  The sheaf $f^{*}(\Omega_{X}^{p-1}\otimes \sL^{p}\otimes \sM)$ is
  ample.  Thus, by Proposition~\ref{ample_subsheaf_of_f*TX} and
  Theorem~\ref{pn_bundle_in_codim_1}, there is an open subset
  $X^\circ\subset X$, with $\codim_X(X\setminus X^\circ)\geq 2$, a
  smooth variety $Y^\circ$, and a $\bP^{d+1}$-bundle
  $\pi^\circ:X^\circ\to Y^\circ$ such that any rational curve from $H$
  meeting $X^\circ$ is a line on a fiber of $\pi^\circ$. Moreover, the
  restriction of $s$ to $X^\circ$ lies in
  $H^0(X^\circ,\wedge^{p}T_{X^\circ/Y^\circ}\otimes
  \sL|_{X^\circ}^{-p}\otimes \sM|_{X^\circ}^{-1})$, and its
  restriction to a general fiber $F$ yields a non-zero section in
  $H^0(F,\wedge^{p}T_{F}\otimes \sL|_{F}^{-p}\otimes \sM|_{F}^{-1})$.
  On the other hand, by Bott's formula, $H^0(\bP^{d+1},\wedge^{p}
  T_{\bP^{d+1}}(-p-1))=0$ unless $p=d+1$.

  Suppose $\dim(Y^\circ)>0$. Since $\codim_X(X\setminus X^\circ)\geq
  2$, $Y^\circ$ contains a complete curve through a general point. Let
  $g:B\to Y^\circ$ be the normalization of a complete curve passing
  through a general point of $Y^\circ$.  Set $X_{B}:=X^\circ
  \times_{Y^\circ} B$, and denote by $\sL_{X_B}$ and $\sM_{X_B}$ the
  pullbacks of $\sL$ and $\sM$ to $X_B$ respectively.  Then $X_{B}\to B$
  is a $\bP^{p}$-bundle, and the section $s$ induces a non-zero section
  in $H^0(X_B,\wedge^{p}T_{X_{B}/B}\otimes \sL_{X_B}^{-p}\otimes
  \sM_{X_B}^{-1})$.  But this is impossible by Corollary~\ref{pbundle}.
  Thus $\dim(Y^\circ)=0$ and $X\simeq \bP^{p}$.
\end{proof}

\begin{cor}\label{vanishing}
  Let $X$ be a smooth projective variety and $\sL$ an ample line bundle
  on $X$.  If $H^0(X,\wedge^{p}T_{X}\otimes \sL^{-p-1-k})\neq 0$ for
  integers $p\ge 1$ and $k\ge 0$, then $k=0$ and $(X,\sL)\simeq
  (\bP^{p},\sO_{\bP^{p}}(1))$.
\end{cor}

\begin{proof}
  Note that $X$ is uniruled by \cite{miyaoka}.  The result follows
  easily from Proposition~\ref{generalvanishing}.
\end{proof}

Here is how we are going to apply these results under the assumptions
of Theorem~\ref{main_thm}.  Suppose that $H^0(X,\wedge^{p}T_{X}\otimes
\sL^{-p})\neq 0$ for some ample line bundle $\sL$ on $X$ and integer
$p\geq 2$.  Then $X$ is uniruled by \cite{miyaoka} and we fix a
minimal covering family $H$ of rational curves on $X$.  Let
$\pi:X^\circ \to Y^\circ$ be the $H$-rational quotient of $X$.  By
shrinking $Y^\circ$ if necessary, we may assume that $Y^\circ$ and
$\pi$ are smooth.  Corollary~\ref{vanishing} provides the vanishing
required to apply Lemma~\ref{2_cases} to $\pi:X^\circ \to Y^\circ$,
yielding the following.

\begin{lemma}\label{fibration}
  Let $Y$ be a smooth variety, $\pi : X\to Y$ a smooth morphism with
  connected fibers, and $\sL$ a line bundle on $X$. Let $F$ be a
  general fiber of $\pi$. Suppose that $F$ is projective and that the
  restriction $\sL|_{F}$ is ample. If $H^0(X,\wedge^{p}T_{X}\otimes
  \sL^{-p})\neq 0$ for some integer $p\geq 2$, then either
  $(F,\sL|_{F})\simeq (\bP^{p-1},\sO_{\bP^{p-1}}(1))$ and
  $H^0(X,\wedge^{p-1}T_{X/Y}\otimes\pi^{*}{T_{Y}}\otimes \sL^{-p})\neq
  0$, or $\dim(F)\geq p$ and $H^0(X,\wedge^{p} T_{X/Y}\otimes
  \sL^{-p})\neq 0$.
\end{lemma}

\begin{proof}
  Corollary~\ref{vanishing} implies that $H^0(F,\wedge^{i}T_F\otimes
  \sL|_F^{-p})= 0$ for $0\leq i\leq p-2$.  So we may apply
  Lemma~\ref{2_cases} with $\sM=\sL^{-p}$ to conclude that either
$H^0(X,\wedge^{p-1}T_{X/Y}\otimes\pi^{*}{T_{Y}}\otimes
    \sL^{-p})\neq 0$, or
$\dim F\geq p$ and $H^0(X,\wedge^{p}T_{X/Y}\otimes \sL^{-p})\neq
    0$.
%
    In the first case we have $H^0(F,\wedge^{p-1}T_{F}\otimes
    \sL|_F^{-p})\neq 0$, and Corollary~\ref{vanishing} implies that
    $(F,\sL|_F)\simeq (\bP^{p-1},\sO_{\bP^{p-1}}(1))$ and so the
    desired statement follows.
\end{proof}

Let $X$, $H$, and $\pi:X^\circ \to Y^\circ$ be as in the above
discussion.  If we are under the first case of Lemma~\ref{fibration},
then Theorem~\ref{pn_bundle_in_codim_1} implies that the
$\bP^{p-1}$-bundle $\pi:X^\circ \to Y^\circ$ can be extended in
codimension $1$.  Next we show that in this case we must have $X\simeq
Q_2$.

\begin{lemma}\label{badbundle}
  Let $X$ be a smooth projective variety and $\sL$ an ample line
  bundle on $X$.  Let $X^\circ\subset X$ be an open subset whose
  complement has codimension at least $2$ in $X$.  Let $\pi :
  X^\circ\to Y^\circ$ be a smooth projective morphism with connected
  fibers onto a smooth quasi-projective variety.  If
  $H^0(X^\circ,\wedge^{p-1}
  T_{X^\circ/Y^\circ}\otimes\pi^{*}{T_{Y^\circ}}\otimes \sL|_{X^\circ}^{-p})\neq
  0$ for some integer $p\geq 2$, then $p=2$, $X^\circ=X\simeq Q_2$,
  and $Y^\circ\simeq \bP^1$.
\end{lemma}

\begin{proof}
  Suppose that for some $p\geq 2$ there is a non-zero section
  $$
  s\in H^0(X^\circ,\wedge^{p-1}
  T_{X^\circ/Y^\circ}\otimes\pi^{*}{T_{Y^\circ}}\otimes \sL|_{X^\circ}^{-p})\neq
  0.
  $$
  By Corollary~\ref{vanishing}, the fibers of $\pi$ are isomorphic to
  $\bP^{p-1}$, and the restriction of $\sL$ to each fiber is
  isomorphic to $\sO_{\bP^{p-1}}(1)$.  Since $\pi$ has relative
  dimension $p-1$, there exists an inclusion $
  \wedge^{p-1}T_{X^\circ/Y^\circ} \otimes\pi^{*}{T_{Y^\circ}}
  \subseteq\wedge^{p}T_{X^\circ}$, and thus $s$, as in (\ref{eq:4}),
  yields a map $ \varphi: \Omega_{X^\circ}^{p-1}\otimes \sL|_{X^\circ}^{p}\to
  T_{X^\circ} $ of rank $p$ at the generic point.  Since
  $\codim_X(X\setminus X^\circ)\geq 2$, $s$ extends to a section
  $\tilde s \in H^0(X,\wedge^{p}T_{X}\otimes \sL^{-p})$. Denote by
  $$
  \widetilde\varphi: \Omega_{X}^{p-1}\otimes \sL^{p}\to T_{X}
  $$
  the associated map, which has rank $p$ at the generic point.

  Let $\sE=\pi_{*}\sL$. By \cite[Corollary 5.4]{fujita75},
  $X^\circ\simeq \bP(\sE)$ over $Y^\circ$ and then
  $\wedge^{p-1}T_{X^\circ/Y^\circ} \otimes \sL^{-p}\simeq \pi^{*}(\det
  \sE^{\dual})$, and $s$ is the pullback of a global section
  $s_{Y^\circ}\in H^0(Y^\circ,T_{Y^\circ}\otimes\det \sE^{\dual})$.
  This implies that the distribution $\sD$ defined by $s$ is
  integrable. Moreover, its leaves are the pullbacks of the leaves of
  the foliation $\sF^\circ$ defined by the map $\det \sE\into
  T_{Y^\circ}$ associated to $s_{Y^\circ}$.

  Since $\codim_X(X\setminus X^\circ)\geq 2$, we can find complete
  curves sweeping out a dense open subset of $Y^\circ$. Let $C$ be a
  general complete curve on $Y^\circ$.  Compactify $Y^\circ$ to a
  smooth variety $Y$, and let $\sF$ be an invertible subsheaf of $T_Y$
  extending $\sF^\circ$.  Then $\sF|_C=\det \sE|_C$ is ample.  By
  \cite[Theorem 0.1]{bogomolov_mcquillan} (see also \cite[Theorem
  1]{KSCT}), the leaf of the foliation $\sF$ through any point of
  $C$ is rational.  We conclude that the leaves of $\sF^\circ$ are
  (possibly noncomplete) rational curves.  Thus the closures of the
  leaves of the distribution $\widetilde \sD$ defined by $\widetilde
  \varphi$ are algebraic.

  Let $F\subset X$ be the closure of a leaf of $\widetilde \sD$ that
  meets $X^\circ$ and let $\eta:\widetilde F\to F$ be its
  normalization.  Then there exists a morphism $\widetilde F\to B$
  onto a smooth rational curve.  The general fiber of this morphism is
  isomorphic to $\bP^{p-1}$ and the restriction of $\sL$ to the
  general fiber is isomorphic to $\sO_{\bP^{p-1}}(1)$.  The fibers are
  thus generically reduced and finally reduced since fibers satisfy
  Serre's condition $S_{1}$.  By \cite[Corollary 5.4]{fujita75},
  $\widetilde F\to B$ is a $\bP^{p-1}$-bundle and, in particular,
  $\widetilde F$ is smooth.

  The section $\tilde s \in H^0(X,\wedge^{p}T_{X}\otimes \sL^{-p})$
  defines a non-zero map $\Omega_{X}^{p}\to \sL^{-p}$.  Since $F$ is the
  closure of a leaf of $\widetilde \sD$ and $\sL|_{F}$ is torsion free, the
  restriction of this map to $F$ factors through a map
  $\Omega_{F}^{p}\to {\sL}|_{F}^{-p}$.  By Lemma
  \ref{extensionpderivation}, this map extends to a map
  $\Omega_{\widetilde F}^{p}\to \eta^{*}{\sL}|_{F}^{-p}$.  Corollary~\ref{pbundle}
  then implies that $p=2$ and $\widetilde F\simeq Q_{2}$.  Moreover
  $\eta^{*}{\sL}|_{F}\simeq \sO_{Q_{2}}(1)$.
  In particular, $\pi : X^\circ\to Y^\circ$ is a $\bP^{1}$-bundle.
  Denote by $H$ the unsplit covering family of rational on $X$ whose 
  general member corresponds to a fiber of $\pi$.

  We claim that the general leaf of $\sF^\circ$ is a complete rational
  curve.  From this it follows that the general leaf of $\widetilde
  \sD$ is compact, and contained in $X^\circ$.  Let $\widetilde F$
  denote the normalization of the closure of a general leaf of
  $\widetilde \sD$.  Since $\widetilde F\simeq Q_{2}$ and
  $\eta^{*}{\sL}|_{F}\simeq \sO_{Q_{2}}(1)$, $X$ admits an unsplit
  covering family $H'$ of rational curves whose general member
  corresponds to a ruling of $\widetilde F\simeq Q_{2}$ that is not
  contracted by $\pi$.  Since $\codim(X\setminus X^\circ)\geq 2$, the
  general member of $H'$ corresponds to a complete rational curve contained in
  $X^\circ$. Its image in $Y^\circ$ is a complete leaf of $\sF^\circ$.
  As we noted above, this implies that $F=\widetilde F\simeq Q_{2}$.  Notice that the section
  $\tilde s$ does not vanish anywhere on a general leaf $F\simeq
  Q_{2}$ of $\sF^\circ$.

  Let $\varphi:X' \to Z'$ be the $(H,H')$-rationally connected quotient of $X$.
  Then the general fiber of $\varphi$ is a leaf $F\simeq Q_{2}$ of
  $\sF^\circ$.  By Lemma~\ref{extending_in_codim_1}, we may assume
  that $\codim_X(X\setminus X')\geq 2$, $Z'$ is smooth, and $\varphi$
  is a proper surjective equidimensional morphism with irreducible and
  reduced fibers.  Therefore $\varphi:X' \to Z'$ is a quadric bundle
  by \cite[Corollary 5.5]{fujita75}. Since the families $H$ and $H'$ are 
  distinct, $\varphi$ is in fact a smooth quadric bundle.

  We claim that in fact $X=F$ and $Z'$ is a point. Suppose otherwise,
  and let $g:C\to Z'$ be the normalization of a complete curve passing
  through a general point of $Z'$. Set $X_C=X' \times_{Z'} C$, denote
  by $\varphi_C:X_C\to C$ the corresponding (smooth) quadric bundle,
  and write $\sL_{X_C}$ for the pullback of $\sL$ to $X_C$.  The
  section $\tilde s$ induces a non-zero section in $H^0(X_C,
  \omega_{X_C/C}^{-1}\otimes \sL^{-2}_{X_C})$ that does not vanish
  anywhere on a general fiber of $\pi_C$.  Thus $ \omega_{X_C/C}^{-1}$
  is ample, contradicting Proposition~\ref{-KX/Y_not_ample}.
\end{proof}



\section{Proof of Theorem {\ref{main_thm}}}
\label{section_proof} 

In order to prove the main theorem, we shall reduce it to the case
when $X$ has Picard number $\rho(X)=1$.  To treat that case, we will
recall some facts about slopes of vector bundles that will be used
later.

\begin{defn}
  Let $X$ be an $n$-dimensional projective variety and $\sH$ an ample
  line bundle on $X$.  Let $\sE$ be a torsion-free sheaf on $X$.  We
  define the slope of $\sE$ with respect to $\sH$ to be
  $\mu_{_\sH}(\sE)=\frac{c_1(\sE)\cdot c_1(\sH)^{n-1}}{\rk(\sE)}$.  We
  say that a torsion-free sheaf $\sF$ on $X$ is
  \emph{$\mu_{_\sH}$-semistable} if for any 
  subsheaf $\sE$ of $\sF$ we have $\mu_{_\sH}(\sE)\leq
  \mu_{_\sH}(\sF)$.
  Given a torsion-free sheaf $\sF$ on $X$, there exists a filtration
  of $\sF$ by (torsion-free) subsheaves
  $$
  0=\sE_0\subsetneq \sE_1\subsetneq \ldots\subsetneq \sE_k=\sF,
  $$
  with $\mu_{_\sH}$-semistable quotients $\cQ_i=\sE_i/\sE_{i-1}$, and
  such that $\mu_{_\sH}(\cQ_1) > \mu_{_\sH}(\cQ_2) > \ldots >
  \mu_{_\sH}(\cQ_k)$.  This is called the \emph{Harder-Narasimhan
    filtration} of $\sF$ \cite{HN75}, \cite[1.3.4]{HuyLehn}.
\end{defn}

\begin{lemma}\label{lem:tf-subsheaf}
  Let $X$ be a smooth $n$-dimensional projective variety and $\sH$ an
  ample line bundle on $X$.  Let $\sF$ be a vector bundle on $X$, $p$
  a positive integer, and $\sN$ an invertible subsheaf of
  $\sF^{\otimes p}$.  Then $\sF$ contains a (torsion-free) subsheaf
  $\sE$ such that $\mu_{_\sH}(\sE) \geq \frac{\mu_{_\sH}(\sN)}{p}$.
\end{lemma}

\begin{proof}
  Consider the Harder-Narasimhan filtration of $\sF$:
  $$
  0=\sE_0\subsetneq \sE_1\subsetneq \ldots\subsetneq \sE_r=\sF,
  $$
  with $\cQ_i=\sE_i/\sE_{i-1}$ $\mu_{_\sH}$-semistable for $1\leq
  i\leq r$, and $\mu_{_\sH}(\cQ_1) > \mu_{_\sH}(\cQ_2) > \ldots >
  \mu_{_\sH}(\cQ_k)$.
  We claim that $\sE=\sE_1=\cQ_1$ satisfies the desired condition.
  In order to prove this, first let $m\in \bN$ be such that
  $\sH^{\otimes m}$ is very ample and let $C\subset X$ be a curve that
  is the intersection of the zero sets of $n-1$ general sections of
  $\sH^{\otimes m}$. Observe that for this curve $C$, and for any
  torsion-free sheaf $\sE$ on $X$,
  \begin{equation}
    \label{eq:6}
    \mu_{_{\sH}}(\sE\resto{C})= m^{n-1}\cdot\mu_{_\sH}(\sE).    
  \end{equation}
  Notice that by abuse of notation we denote the restriction of $\sH$
  to $C$ by the same symbol. Let $\sG_i=\sE_i\resto C$ and
  $\cP_i=\cQ_i\resto C$.  By the Mehta-Ramanathan Theorem
  (\cite[6.1]{MR_semistable_curves}, \cite[7.2.1]{HuyLehn}) the
  Harder-Narasimhan filtration of $\sF\resto{C}$ is exactly the
  restriction to $C$ of the Harder-Narasimhan filtration of $\sF$ (we
  may assume that $m$ was already chosen large enough for this theorem
  to apply as well):
  $$
  0=\sG_0\subsetneq \sG_1\subsetneq \ldots\subsetneq \sG_r=\sF\resto C.
  $$
  As $X$ is smooth, so is $C$ and hence all torsion-free sheaves on
  $C$, in particular the $\sG_i$ and the $\cP_i$, are locally free.
  Then for each $1\leq i\leq r$ there exists a filtration
  $$
  \sG_{i}^{\otimes p}=\sG_{i, 0}\supseteq \sG_{i,1}\supseteq
  \ldots\supseteq \sG_{i, p}\supseteq \sG_{i, p+1}=0,
  $$
  with quotients $\sG_{i,j}/\sG_{i,j+1}\simeq \sG_{i-1}^{\otimes
    j}\otimes \cP_{i}^{\otimes (p-j)}$.
  From these filtrations, we see that the inclusion $\sN\into
  \sF^{\otimes p}$ induces an inclusion $\sN\resto C\into
  \cP_1^{\otimes i_1}\otimes \ldots \otimes \cP_{k}^{\otimes i_k}$,
  for suitable non-negative integers $i_j$'s such that $\sum i_j = p$.
  Since each $\cP_i$ is $\mu_{_\sH}$-semistable (on $C$), so is the
  tensor product $\cP_1^{\otimes i_1}\otimes \ldots \otimes
  \cP_{k}^{\otimes i_k}$ \cite[Theorem 3.1.4]{HuyLehn}.  Hence
  $$
  \mu_{_\sH}(\sN) = \frac{\mu_{_\sH}(\sN\resto C)}{m^{n-1}} \leq
  \frac{\mu_{_\sH}(\cP_1^{\otimes i_1}\otimes \ldots \otimes
    \cP_{k}^{\otimes i_k})}{m^{n-1}} = \frac{\sum i_j
    \mu_{_\sH}(\cP_j)}{m^{n-1}} \leq \frac{p\mu_{_\sH}(\cQ_1\resto
    C)}{m^{n-1}}= p\mu_{_\sH}(\cQ_1),
  $$
  and so $\sE=\sE_1=\cQ_1$ does indeed satisfy the required property.
\end{proof}

Now we can prove our main theorems.

\begin{thm}\label{beauville_for_rho=1}
  Let $X$ be a smooth $n$-dimensional projective variety with
  $\rho(X)=1$, $\sL$ an ample line bundle on $X$, and $p$ a positive
  integer.  Suppose that $H^0(X,T_{X}^{\otimes p}\otimes \sL^{-p})\neq
  0$. Then either $(X,\sL)\simeq (\p^{n},\sO_{\p^{n}}(1))$, or $p=n\geq
  3$ and $(X,\sL)\simeq(Q_{p},\sO_{Q_{p}}(1))$.
\end{thm}

\begin{proof}
  First notice that $X$ is uniruled by \cite{miyaoka}, and hence a
  Fano manifold with $\rho(X)=1$. The result is clear if $\dim X=1$,
  so we assume that $n\geq 2$.  Fix a minimal covering family $H$ of
  rational curves on $X$.
  By Lemma~\ref{lem:tf-subsheaf}, $T_X$ contains a torsion-free
  subsheaf $\sE$ such that $\mu_{_\sL}(\sE) \geq
  \frac{\mu_{_\sL}(\sL^p)}{p} = \mu_{_\sL}(\sL)$.  This implies that
  $\frac{\deg f^*\sE}{\rk \sE}\geq \deg f^*\sL$ for a general member
  $[f]\in H$.  If $r=\rk(\sE)=1$, then $\sE$ is ample and we are done
  by Wahl's theorem.  Otherwise, as $f^*\sE$ is a subsheaf of
  $f^*T_X\simeq \sO_{\bP^1}(2)\oplus \sO_{\bP^1}(1)^{\oplus d}\oplus
  \sO_{\bP^1}^{\oplus (n-d-1)}$, we must have $\deg f^*\sL=1$ and
  either $f^*\sE$ is ample, or $f^*\sE\simeq \sO_{\bP^1}(2)\oplus
  \sO_{\bP^1}(1)^{\oplus r-2}\oplus \sO_{\bP^1}$ for a general $[f]\in
  H$.
  If $f^*\sE$ is ample, then $X\simeq \p^n$ by
  Proposition~\ref{ample_subsheaf_of_f*TX}, using the fact that
  $\rho(X)=1$.  If $f^*\sE$ is not ample, then $\sO_{\bP^1}(2)\subset
  f^*\sE$ for general $[f]\in H$, and so $\cC_x\subset
  \p(\sE^*\otimes\kappa(x))$ for a general $x\in X$. Thus by
  \cite[2.6]{artigo_tese} $(f^*T_X^+)_o\subset (f^*\sE)_o$ for a
  general $o\in \p^1$ and a general $[f]\in H$.  Since $f^*T_X^+$ is a
  subbundle of $f^*T_X$, we have an inclusion of sheaves
  $f^*T_X^+\into f^*\sE$, and thus $\det(f^*\sE)=f^*\omega_X^{-1}$.
  Since $\rho(X)=1$, this implies that $\det \sE^{**}=\omega_X^{-1}$,
  and thus $0\neq h^{0}(X,\wedge^{r} T_{X}\otimes
  \omega_{X})=h^{n-r}(X,\sO_{X})$.  The latter is zero unless $n=r$
  since $X$ is a Fano manifold.  If $n=r$, then we must have
  $\omega_X^{-1}\simeq \sL^{\otimes n}$.  Hence $X\simeq Q_n$ by
  \cite{kobayashi_ochiai}.
\end{proof}

\begin{proof}[{Proof of Theorem {\ref{main_thm}}}]
  Let $X$ be a smooth projective variety and $\sL$ an ample line
  bundle on $X$ such that $H^0(X,\wedge^{p} T_{X}\otimes \sL^{-p})\neq
  0$.  By Theorem~\ref{beauville_for_rho=1}, we may assume that
  $\rho(X)\geq 2$.  We may also assume that $p\geq 2$ as the case
  $p=1$ is just Wahl's theorem.  We shall proceed by induction on $n$.

  Notice that $X$ is uniruled by \cite{miyaoka}.  Let
  $H\subset\rat(X)$ be a minimal covering family of rational curves on
  $X$, and $[f]\in H$ a general member.  By analyzing the degree of
  the vector bundle $f^*(\wedge^{p}T_{X}\otimes \sL^{-p})$, we
  conclude that $f^*\sL\simeq \sO_{\p^1}(1)$, and thus $H$ is unsplit.
  Let $\pi^\circ:X^\circ \to Y^\circ$ be the $H$-rational quotient of
  $X$.  By shrinking $Y^\circ$ if necessary, we may assume that
  $\pi^\circ$ is smooth.  Since $\rho(X)\geq 2$, we must have $\dim
  Y^\circ\geq 1$ by \cite[IV.3.13.3]{kollar}.

  Let $F$ be a general fiber of $\pi^\circ$ and set $k=\dim F$.
  By Lemma~\ref{fibration}, either
  \begin{itemize}
  \item $k=p-1$, $(F,\sL|_{F})\simeq (\p^{p-1},\sO_{\p^{p-1}}(1))$,
    and $H^0(X^\circ,\wedge^{p-1}T_{X^\circ/Y^\circ}\otimes
    \pi^{*}{T_{Y^\circ}}\otimes \sL^{-p})\neq 0$, or
  \item $k\geq p$ and $H^0(X^\circ,\wedge^{p}
    T_{X^\circ/Y^\circ}\otimes \sL^{-p})\neq 0$.
   \end{itemize}

   In the first case $\pi:X^\circ \to Y^\circ$ is a $\p^{p-1}$-bundle
   and we may assume that $\codim_X(X\setminus X^\circ)\geq 2$
   by Theorem~\ref{pn_bundle_in_codim_1}.  Then we apply
   Lemma~\ref{badbundle} and conclude that $X\simeq Q_2$.

   In the second case, the induction hypothesis implies that either
   $(F,\sL|_{F})\simeq (\p^{k},\sO_{\p^{k}}(1))$, or $k=p$ and
   $(F,\sL|_{F})\simeq (Q_{p},\sO_{Q_{p}}(1))$.  If $F\simeq \p^k$,
   again by Theorem~\ref{pn_bundle_in_codim_1}, $\pi:X^\circ \to
   Y^\circ$ is a $\p^{k}$-bundle, and we may assume that
   $\codim_X(X\setminus X^\circ)\geq 2$.  As in the end of the proof
   of Proposition~\ref{generalvanishing}, we reach a contradiction by
   applying Corollary~\ref{pbundle} to $X^\circ \times_{Y^\circ} B\to B$,
   where $B\to Y^\circ$ is the normalization of a complete curve
   passing through a general point of $Y^\circ$.

   Suppose now that $F\simeq Q_p$. Then, by
   Lemma~\ref{extending_in_codim_1} and \cite[Corollary
   5.5]{fujita75}, $\pi^\circ$ can be extended to a quadric bundle
   $\pi:X' \to Y'$ with irreducible and reduced fibers,
   where $X'$ is an open subset of $X$ with  
   $\codim_X(X\setminus X')\geq 2$, and $Y'$ is smooth.
   Denote by $X''$ the open
   subset of $X'$ where $\pi$ is smooth.  Notice that
   $\codim_{X'}(X'\setminus X'')\geq 2$.  A non-zero global section of
   $\wedge^{p}T_{X}\otimes \sL^{-p}$ restricts to a non-zero global
   section of $\wedge^{p}T_{X''/Y'}\otimes \sL|_{X''}^{-p}$, which, in
   turn, extends to a non-zero global section $s\in H^0(X',
   \omega_{X'/Y'}^{-1}\otimes \sL|_{X'}^{-p})$ since $X'$ is smooth.  The
   section $s$ does not vanish anywhere on a general fiber of $\pi$.
   
   Let $g:C\to Y'$ be the normalization of a complete curve passing
   through a general point of $Y'$. Set $X_C=X' \times_{Y'} C$, denote
   by $\pi_C:X_C\to C$ the corresponding quadric bundle, and write
   $\sL_{X_C}$ for the pullback of $\sL$ to $X_C$.  The general fiber
   of $\pi_C$ is smooth. Now notice that $X_C$ is a local complete
   intersection variety, and nonsingular in codimension one, since the
   fibers of $\pi$ are reduced.  In particular, $X_C$ is a normal
   Gorenstein variety, and the morphism $\pi_C$ is generically smooth.
   The section $s$ induces a non-zero section in $H^0(X_C,
   \omega_{X_C/C}^{-1}\otimes \sL^{-p}_{X_C})$ that does not vanish
   anywhere on the general fiber of $\pi_C$.  Thus $
   \omega_{X_C/C}^{-1}$ is ample, contradicting
   Proposition~\ref{-KX/Y_not_ample}.
\end{proof}


\providecommand{\bysame}{\leavevmode\hbox to3em{\hrulefill}\thinspace}
\providecommand{\MR}{\relax\ifhmode\unskip\space\fi MR}
\providecommand{\MRhref}[2]{%
  \href{http://www.ams.org/mathscinet-getitem?mr=#1}{#2}
}
\providecommand{\href}[2]{#2}

\end{document}